\documentclass[a4paper, 11pt]{amsart}
\usepackage[latin1]{inputenc}
\usepackage[T1]{fontenc}
\usepackage[british]{babel}
\usepackage{amssymb}
\usepackage{amsmath}
\usepackage{amsthm}
\usepackage{amscd}
\usepackage{amsfonts}
\usepackage{stmaryrd}
\usepackage{pb-diagram}
\usepackage{epic,eepic,epsfig}
\usepackage{a4wide}
\usepackage{nextpage}
\usepackage{fancyhdr}

\newtheorem{prop}{Proposition}[section]
\newtheorem{cor}[prop]{Corollary}

\newtheorem{lem}[prop]{Lemma}

\newtheorem{conj}[prop]{Conjecture}

\newcommand{\Jac} {\mathop{\mathrm{Jac}}}

\newcommand{\Sym} {\mathop{\mathrm{Sym}}}

\newcommand{\degr} {\mathop{\mathrm{deg}}}

\newcommand{\Isog} {\mathop{\mathrm{Isog}}}
\newcommand{\Prep} {\mathop{\mathrm{Prep}}}

\begin{document}

\title[]{Zhang's conjecture and squares of abelian surfaces}


\author[Fabien {\sc Pazuki}]{{\sc Fabien} Pazuki}
\address{Fabien {\sc Pazuki}\\
IMB U. Bordeaux 1\\
351, cours de la Lib\'eration\\
33405 Talence, France\\}
\email{fabien.pazuki@math.u-bordeaux1.fr}
\urladdr{http://www.math.u-bordeaux1.fr/~pazuki}

\maketitle

\begin{center}
 January 26th, 2010
\end{center}

\begin{centering}
{\small{\textsc{R\'esum\'e}: On donne dans cette note des exemples de carr\'es de surfaces ab\'eliennes violant la conclusion de la conjecture de Zhang sur l'intersection des sous-vari\'et\'es et des points pr\'ep\'eriodiques.}}
\end{centering}

\vspace{0.3cm}

\begin{centering}
{\small{\textsc{Abstract}: We give in this paper some squares of abelian surfaces that are counterexamples to a conjecture formulated by Zhang about the intersection of subvarieties and preperiodic points.}}
\end{centering}

{\flushleft
}

\thispagestyle{empty}

\section{Introduction}
D. Ghioca and T. Tucker found a family of counterexamples to Zhang's dynamical Manin-Mumford Conjecture 1.2.1 of \cite{Zha}. They use squares of elliptic curves with complex multiplication. D. Ghioca asked whether this counterexample could be generalized. We present here counterexamples of greater dimension. We recall a few definitions: an endomorphism $\varphi:X\to X$ of a projective variety is said to have a \textit{polarization} if there exists an ample divisor $D$ such that $\varphi^{*}D\sim d D$ for some $d>1$, where $\sim$ stands for the linear equivalence. A subvariety $Y$ of $X$ is \textit{preperiodic} under $\varphi$ if there exists integers $m\geq0$ and $k>0$ such that $\varphi^{m+k}(Y)=\varphi^{m}(Y)$. We denote $\Prep_{\varphi}(X)$ the set of preperiodic points of $X$ under the action of $\varphi$. We now recall the conjecture:

\begin{conj}\label{Zhang}
(Algebraic Dynamical Manin-Mumford) Let $\varphi :X\to X$ be an endomorphism of a projective variety defined over a number field $K$ with a polarization, and let $Y$ be a subvariety of $X$. If $Y\cap \Prep_{\varphi}(X)$ is Zariski-dense in $Y$, then $Y$ is a preperiodic subvariety.
\end{conj}

We will use the following lemma of Ghioca and Tucker \cite{GhTuZh}:
\begin{lem}\label{ghioca tucker}
Let $A$ be a simple abelian variety and $\varphi_{1},\varphi_{2}\in{\Isog(A)}$ be nonzero. Let $m\geq0$ and $k>0$ be two integers. Let $\Delta=\{(x,x)|x\in{A}\}$. Then $(\varphi_{1}^m,\varphi_{2}^m)(\Delta)\subset (\varphi_{1}^{m+k},\varphi_{2}^{m+k})(\Delta)$ if and only if $\varphi_{1}^k=\varphi_{2}^k$.
\end{lem}
In a nutshell, for the examples that we find in this paper, $\Delta$ is preperiodic under $\varphi=(\varphi_{1},\varphi_{2})$ if and only if the action of $\varphi_{1}$ and $\varphi_{2}$ differ by a root of unity.
\begin{proof}
For the direct part, one checks that $\varphi_{1}^{k}-\varphi_{2}^{k}$ sends nontorsion points to torsion points, hence is the zero map. The converse comes from the surjectivity of $\varphi_{1}$ and $\varphi_{2}$.
\end{proof}

The examples we provide are inspired by Ghioca and Tucker's original ones, the main problem to overcome is finding a polarizable situation. The idea used here is the theorem of the cube on abelian varieties combined with some particular properties of the field of definition.

For other reflections on the dynamical Manin-Mumford conjecture, one can refer to \cite{GhTuZh} or \cite{YuZh}.

\begin{centering}
{{\textbf{Aknowledgements}: Many thanks to D. Ghioca for fruitful conversations and for the pleasant stay I had at the University of Lethbridge. I also thank the anonymous referee for his helpful comments.}}
\end{centering}

\section{Polarizability criterion}

We give in this section a few formulas useful to get information on the weight of complex multiplication. We start with a general fact:

\begin{prop}\label{cube}
Let $A$ be an abelian variety, $V$ a variety and $f,g,h$ three morphisms from $V$ to $A$. Then for any divisor $D\in{Div(A)}$, one has
$$(f+g+h)^{*}D-(f+g)^{*}D-(g+h)^{*}D-(f+h)^{*}D+f^{*}D+g^{*}D+h^{*}D\sim 0.$$
\end{prop}
\begin{proof}
This statement is a direct consequence of the theorem of the cube. For a proof, see for example \cite{HiSi}, Corollary A.7.2.4 page 123.
\end{proof}

Let $A$ be an abelian variety and suppose it has complex multiplication by a ring $R$, \textit{i.e.} the ring of endomorphisms of $A$ contains $R$ and $R$ contains stricly $\mathbb{Z}$ (see \cite{Sil}, chapter II.1 for the case of elliptic curves). Then we have the following lemma:

\begin{lem}
Let $A$ be an abelian variety and let $D$ be a divisor on $A$. Let $n\in\mathbb{N}$ and $\alpha\in{R}$. Then
\begin{equation}
[n+\alpha]^{*}D\sim n[1+\alpha]^{*}D-(n-1)[\alpha]^{*}D+n(n-1)D.
\end{equation}
\end{lem}
\begin{proof}
Use Proposition~\ref{cube} with $f=[n-1]$, $g=[\alpha]$ and $h=[1]$. The result follows by a recurrence and a telescoping sum.
\end{proof}

\begin{cor}\label{n+alpha}
If one chooses $D$ such that $[\alpha]^{*}D\sim D$, then $[n+\alpha]$ is polarized by $D$ if and only if $[1+\alpha]$ is polarized by $D$ and one has 
\begin{equation}\label{n plus alpha}
[n+\alpha]^{*}D\sim n[1+\alpha]^{*}D+(n-1)^2 D.
\end{equation}
\end{cor}

\section{Theta divisor in dimension 2}

Let $C$ be a curve of genus 2 defined over $\bar{\mathbb{Q}}$. Choose an affine equation $y^{2}=f(x)$ with $\degr(f)=5$, and let $\infty$ be the point at infinity. Let $\Jac(C)$ denote the jacobian of $C$. We denote by $cl(D)$ the linear equivalence class of any divisor $D$ on $\Jac(C)$. Let $\Theta=j(C)$ the theta divisor, where
\begin{align*}
j :C & \hookrightarrow \Jac(C)\\
P & \longrightarrow cl((P)-(\infty)).
\end{align*}

Consider the surjective map
\begin{align*}
\Sym\,\!^{2}(C) & \longrightarrow \Jac(C)\\
\{P_{1},P_{2}\}& \longrightarrow cl((P_{1})+(P_{2})-2(\infty)).
\end{align*}

Take $P_{1}=(x_{1},y_{1})$ and $P_{2}=(x_{2},y_{2})$. Then the hyperelliptic involution $\iota:(x,y)\to (x,-y)$ gives the multiplication by $[-1]$ on $\Jac(C)$. We have to blow down all the points $\{P,\iota(P)\}$ on $\Sym\,\!^{2}(C)$ to the origin of $\Jac(C)$. The theta divisor is then the image of the set of all pairs $\{P,\infty\}$ with $P\in{C}$. It is symmetric and ample.

If one chooses the affine equation to be $y^2=f(x)$ with $\degr(f)=6$, then there are two points at infinity $\infty^{+}$ and $\infty^{-}$, and over a field extension one gets a point $\infty$ such that $(\infty^{+})+(\infty^{-})\sim 2(\infty)$. Then one would work with the divisor $D_{0}=t_{\infty-\infty^{+}}^{*}\Theta+t_{\infty-\infty^{-}}^{*}\Theta\sim 2\Theta$, where $t_{P}$ stands for the translation by the point $P$.

\section{Example in degree $5$}

Let us focus on the curve $C$ with affine model $y^2=x^5-x$. In this particular case, we get a jacobian with complex multiplication, coming from $[i]:(x,y)\rightarrow (-x,iy)$, where $i^2=-1$.

We have $[i]^{*}(x,y)=(-x,-iy)$ on the curve, which gives $[i]^{*}(\{(x,y),\infty\})=\{(-x,-iy),\infty\}$ on $\Sym\,\!^{2}(C)$, thus $[i]^{*}\Theta\sim\Theta$. 
Let us use Proposition~\ref{cube} in the following situation: $A=V=\Jac(C)$, $D=\Theta$, $f=[i]$, $g=[1]$ and $h=[-1]$. Then we get $$[i]^{*}\Theta-[1+i]^{*}\Theta-[i-1]^{*}\Theta+[i]^{*}\Theta+\Theta+[-1]^{*}\Theta\sim0,$$
thus using $[-1]^{*}\Theta \sim \Theta$ and $[i]^{*}\Theta \sim \Theta$ we have
\begin{equation}\label{eq}
[1+i]^{*}\Theta+[1-i]^{*}\Theta\sim 4\Theta.
\end{equation}

Let us now remark that $i(1-i)=1+i$, so $[1+i]^{*}\Theta\sim[i(1-i)]^{*}\Theta\sim[1-i]^{*}[i]^{*}\Theta\sim[1-i]^{*}\Theta.$ Using this in equation (\ref{eq}), one gets

\begin{equation}\label{1+i}
[1+i]^{*}\Theta\sim 2\Theta.
\end{equation}

Then using Corollary~\ref{n+alpha} one gets $[2+i]^{*}\Theta\sim 5\Theta$ and $[2-i]^{*}\Theta\sim 5\Theta$.

Let $A=\Jac(C)$ and let $\varphi=[2+i]\times [2-i]$. Consider the following situation
$$\varphi : A\times A\longrightarrow A\times A.$$
The morphism $\varphi$ is polarized by $D=\pi_{1}^{*}\Theta+\pi_{2}^{*}\Theta$, where $\pi_{1}$ and $\pi_{2}$ are respectively the first and second projections, and $\varphi^{*}D\sim 5 D$. Choose $Y=\Delta=\{(P,P)\in{A\times A}\}$ to be the diagonal. The intersection $Y\cap\mathrm{Prep(A\times A)}$ is Zariski-dense in $Y$. Then $\varphi^k Y=Y$ implies that for every $P\in{Y}$ we have $[2-i]^k P=[2+i]^k P$. But $\frac{2-i}{2+i}$ is not a root of unity. Use lemma \ref{ghioca tucker}. We thus have provided a square of an abelian surface that contradicts Conjecture \ref{Zhang}.

\section{Example in degree $6$}

Let us focus on the curve $C$ with affine model $y^2=x^6-1$. With this choice of affine model, one has two points at infinity denoted $\infty^{+}$ and $\infty^{-}$. We consider the endomorphism $[\alpha]:(x,y)\rightarrow (\alpha x,y)$, where $\alpha^6=1$. This morphism gives rise to a complex multiplication endomorphism on the surface $\Jac(C)$ that will also be denoted $[\alpha]$. The divisor $(\infty^{+})+(\infty^{-})$ is invariant under $[\alpha]$. We split the study into two cases, whether we have $\alpha=j$ where $j^2+j+1=0$ or $\alpha^2-\alpha+1=0$. We begin with $\alpha=j$. Define the divisor $$D_{1}=D_{0}+[j]^{*}D_{0}+[j^2]^{*}D_{0}+[-1]^{*}D_{0}+[-j]^{*}D_{0}+[-j^2]^{*}D_{0}.$$ One verifies that $[\pm j^{m}]^{*}D_{1}\sim D_{1}$ for $m=0,1,2$. Let us use Proposition \ref{cube} in the following situation: $A=V=\Jac(C)$, $D=D_{1}$, $f=[1]$, $g=[j]$ and $h=[j]$. Then we get $$[1+2j]^{*}D_{1}-2[1+j]^{*}D_{1}-[2j]^{*}D_{1}+D_{1}+2[j]^{*}D_{1}\sim0,$$
thus using $1+2j=j-j^2$ and $1+j=-j^2$, plus $[-1]^{*}D_{1} \sim D_{1}$ and $[j]^{*}D_{1} \sim D_{1}$ we have
\begin{equation}\label{1-j}
[1-j]^{*}D_{1}\sim 3D_{1}.
\end{equation}

Let us now remark that $(1-j^2)j=j-1$, so $[1-j]^{*}D_{1}\sim[j(1-j^2)]^{*}D_{1}\sim[1-j^2]^{*}[j]^{*}D_{1}\sim[1-j^2]^{*}D_{1}.$ Using this in equation (\ref{1-j}), one gets

\begin{equation}\label{1-j^2}
[1-j^2]^{*}D_{1}\sim 3D_{1}.
\end{equation}

Then by using (\ref{1-j}) and Corollary~\ref{n+alpha} one gets $[2-j]^{*}D_{1}\sim 7D_{1}$, and by using (\ref{1-j^2}) and Corollary~\ref{n+alpha} one gets $[2-j^2]^{*}D_{1}\sim 7D_{1}$.

Let $A=\Jac(C)$ and let $\varphi=[2-j]\times [2-j^2]$. Consider the following situation
$$\varphi : A\times A\longrightarrow A\times A.$$
The morphism $\varphi$ is polarized by $D=\pi_{1}^{*}D_{1}+\pi_{2}^{*}D_{1}$, where $\pi_{1}$ and $\pi_{2}$ are respectively the first and second projections, and $\varphi^{*}D\sim 7 D$. Choose $Y=\Delta=\{(P,P)\in{A\times A}\}$ to be the diagonal. The intersection $Y\cap\mathrm{Prep(A\times A)}$ is Zariski dense in $Y$. Then $\varphi^k Y=Y$ implies that for every $P\in{Y}$ we have $[2-j]^k P=[2-j^2]^k P$. But $\frac{2-j}{2-j^2}$ is not a root of unity.

One may deal with the case $\alpha^2=\alpha-1$ in the same way, using Proposition \ref{cube} with $f=[1]$, $g=h=[-\alpha]$.

\section{Multiplication by $\zeta_{5}$ not polarized by $\Theta$}

Let us focus on the curve $C$ with affine model $y^2=x^5-1$. In this particular case, we get a jacobian with complex multiplication coming from $[\zeta_{5}]:(x,y)\rightarrow (\zeta_{5}x,y)$, where $\zeta_{5}^5=1$. We have $[\zeta_{5}]^{*}(x,y)=(\zeta_{5}^4 x,y)$ on the curve, which gives $[\zeta_{5}]^{*}(\{(x,y),\infty\})=\{(\zeta_{5}^4 x,y),\infty\}$ on $\Sym\,\!^{2}(C)$, thus $[\zeta_{5}]^{*}\Theta\sim\Theta$. We gather a few pullback formulas in this particular setting:

\begin{lem}
Let $m$ and $n$ be integers. One has
\begin{equation}
[n+\zeta_{5}^m]^{*}\Theta+[n-\zeta_{5}^m]^{*}\Theta\sim (2n^2+2)\Theta.
\end{equation}
\begin{equation}\label{somme}
[1+\zeta_{5}]^{*}\Theta+[1+\zeta_{5}^2]^{*}\Theta\sim 3\Theta.
\end{equation}
\begin{equation}\label{produit}
[(1+\zeta_{5})(1+\zeta_{5}^2)]^{*}\Theta\sim \Theta.
\end{equation}
\end{lem}
\begin{proof}
The first equality can be deduced from (\ref{n plus alpha}), using $\alpha=\zeta_{5}^m$ and $\alpha=-\zeta_{5}^m$. The second equality comes from the application of Proposition \ref{cube} with $f=[1]$, $g=[\zeta_{5}]$ and $h=[\zeta_{5}^2]$. The last equality comes from the relation $1+\zeta_{5}+\zeta_{5}^2+\zeta_{5}^3=-\zeta_{5}^4$ and the fact that $[\zeta_{5}]^{*}\Theta\sim \Theta$.
\end{proof}
As opposed to the endomorphisms $[i]$ and $[j]$ in the first examples, the formulas (\ref{somme}) and (\ref{produit}) show that $[\zeta_{5}]$ will not be polarized by $\Theta$. 

%
%
%
%
%

\vfill

\end{document}